\documentclass[preprint,11pt]{elsarticle}

\usepackage{amssymb}
\usepackage{amsmath}
\usepackage{amsfonts}
\usepackage{amssymb}
\usepackage{indentfirst,latexsym,bm}
\usepackage{amsthm}
\usepackage{color,xcolor}
\usepackage[all]{xy}
\input{mathrsfs.sty}
\newtheorem{theorem}{Theorem}[section]
\newtheorem{lemma}[theorem]{Lemma}
\newtheorem{corollary}[theorem]{Corollary}

\theoremstyle{definition}
\newtheorem{definition}{Definition}[section]

\theoremstyle{definition}
\newtheorem{example}{Example}[section]

\theoremstyle{remark}
\newtheorem{remark}{Remark}[section]
\theoremstyle{question}

\numberwithin{equation}{section}

\journal{XXX}

\begin{document}

\begin{frontmatter}



\title{The polar decomposition of the product of three operators}
\author[shnu]{Dingyi Du}
\ead{dudingyi234@163.com}
\author[shnu]{Qingxiang Xu}
\ead{qingxiang$\_$xu@126.com}
\author[shnu]{Shuo Zhao}
\ead{z13427525400@163.com}
\address[shnu]{Department of Mathematics, Shanghai Normal University, Shanghai 200234, PR China}

\begin{abstract}
In the setting of adjointable operators on Hilbert $C^*$-modules, this paper deals with the polar decomposition of the product of three operators. The relationship between the polar decompositions associated with three operators is clarified. Based on this relationship, a formula for the polar decomposition of a multiplicative perturbation of an operator is provided. In addition, some characterizations of the polar decomposition associated with three operators are provided.
\end{abstract}


\end{frontmatter}



\section*{Introduction}

Given Hilbert modules $H$ and $K$ over a $C^*$-algebra $\mathfrak{A}$, let $\mathcal{L}(H,K)$ denote the set of adjointable operators from $H$ to $K$, with the abbreviation $\mathcal{L}(H)$ whenever $H=K$. A multiplicative perturbation of an operator $T\in\mathcal{L}(H,K)$ can be represented as $M=ETF^*$, where $E\in\mathcal{L}(K)$ and $F\in\mathcal{L}(H)$ might be not invertible. In the matrix case, norm estimations of $U_M-U_T$ are carried out via the formulas for $U_M$ and $U_T$ derived from the Singular Value Decomposition (SVD) of $M$ and $T$ respectively \cite{Chen-Li-2,Chen-Li-Sun,Hong-Meng-Zheng}, in which $U_M$ (resp.\,$U_T$) is the partial  isometry factor of the polar decomposition $M$ (resp.\,$T$).  As far as we know, little has been done in norm estimations of $U_M-U_T$ when $M$ and $T$ are two general operators acting on an infinite-dimensional space. One of the main difficulties on this issue is the derivation of the formulas for $U_M$, which apparently can not be fulfilled via the SVD, since it works only for matrices. This leads us to study the polar decomposition of the product of three operators.

 For every adjointable operator $B$ on a Hilbert $C^*$-module, let $|B|$ denote the square root of $B^*B$. That is, $|B|=(B^*B)^\frac12$.
 To give characterizations of the centered operator introduced in \cite{MM}, the polar decomposition of the product of two operators is dealt with in \cite[Section~2]{Ito-IEOT} for Hilbert space operators.
Let $T=U_T|T|$ and $S=U_S|S|$ be the polar decompositions. Suppose that the polar decomposition of $|T|\cdot |S^*|$ is given by $|T|\cdot |S^*|=W_1\big||T|\cdot |S^*|\big|$. Then it is proved in \cite[Theorem~2.1]{Ito-IEOT} that $TS=U_TW_1U_S|TS|$ will be the polar decomposition. The reverse  is also considered in \cite[Lemma~3.1]{LLX-BJMA}, which indicates that $|T|\cdot|S^*|=U_T^*W_2U_S^*\big||T|\cdot |S^*|\big|$ will be the polar decomposition whenever $TS=W_2|TS|$ is the polar decomposition. As an application, three equivalent conditions are stated in \cite[Theorem~3.5]{LLX-BJMA}, which is a generalization of \cite[Theorem~2.3]{Ito-IEOT}. Consequently, a generalization of \cite[Theorem~3.1]{Ito-JOT} is obtained in \cite[Theorem~4.14]{LLX-BJMA}, where
a characterization of the binormal operator \cite{Campbell-1}  is provided in terms of the generalized Aluthge transform.

The purpose of this paper is to set up the theory for the polar decomposition of a multiplicative perturbation of an operator, and meanwhile to give some non-trivial generalizations of the results  obtained originally in \cite{Ito-IEOT,Ito-JOT,LLX-BJMA}.

The paper is organized as follows. Some basic knowledge about Hilbert $C^*$-modules and adjointable operators are provided in Section~\ref{sec:preliminaries}. Section~\ref{sec:equivalent characterization} is devoted to the study of the relationship between the polar decompositions associated with three operators.
Let $H$ and $K$ be Hilbert $\mathfrak{A}$-modules, and let  $T\in \mathcal{L}(K), A\in\mathcal{L}(H,K)$ and $S\in\mathcal{L}(H)$ be such that all of them have the polar decompositions. Put
\begin{equation}\label{equ:defn of X and Y}X=TAS\quad\mbox{and}\quad Y=|T|\cdot A\cdot |S^*|.\end{equation}
In Theorem~\ref{thm:main result of 3 operators}, we will show that $X$ has the polar decomposition $X=U_X|X|$ if and only if $Y$ has the polar decomposition $Y=U_Y|Y|$, and in which case the relationship between $U_X$ and $U_Y$ described as above for $A=I_H$   is  also valid for every $A\in\mathcal{L}(H,K)$. Thus, a full generalization of \cite[Lemma~3.1]{LLX-BJMA} is obtained.

Two applications of Theorem~\ref{thm:main result of 3 operators} are dealt with in the remaining parts of this paper.  Based on Theorem~\ref{thm:main result of 3 operators} together with certain block matrix technique, a formula for $U_M$   in the unitarily equivalent sense is derived for the first time in the case that $M$ is a multiplicative perturbation of an operator; see  Theorem~\ref{thm:polar decomposition for multiplicative perturbation} for the details.
Another application of Theorem~\ref{thm:main result of 3 operators} is to find out new characterizations of the polar decomposition associated with three operators. This is fulfilled in Theorem~\ref{thm:technical lemma for three operators}, which serves as the main result of Section~\ref{sec:three operators}. Based on
Theorem~\ref{thm:technical lemma for three operators}, some non-trivial generalizations of \cite[Theorems~3.5 and 4.14]{LLX-BJMA} are obtained; see
Theorem~\ref{thm:application of main thm} and Corollary~\ref{cor:equivalent conditions of binormal} for the details.

\section{Preliminaries}\label{sec:preliminaries}
Hilbert $C^*$-modules are generalizations of Hilbert spaces by allowing inner products to take values in some $C^{*}$-algebras instead of the complex  field \cite{Lance,Manuilov,Paschke}. Given  Hilbert $\mathfrak{A}$-modules $H$ and $K$ over a $C^*$-algebra $\mathfrak{A}$, let $\mathcal{L}(H,K)$
be the set of operators $T:H\to K$ for which there is an operator $T^*:K\to
H$ such that $$\langle Tx,y\rangle=\langle x,T^*y\rangle \quad \mbox{for every $x\in H$ and $y\in K$}.$$  We call ${\mathcal
L}(H,K)$ the set of adjointable operators from $H$ to $K$. For every
$T\in \mathcal{L}(H,K)$, its  range and  null space are denoted by
${\mathcal R}(T)$ and ${\mathcal N}(T)$, respectively. In case
$H=K$, $\mathcal{L}(H,H)$ which is abbreviated to $\mathcal{L}(H)$, is a
$C^*$-algebra. Let $\mathcal{L}(H)_+$ \big(resp.\,$\mathcal{L}(H)_{\mbox{sa}}$\big) denote the set of positive (resp.\,self-adjoint) elements
in $\mathcal{L}(H)$.  When $A\in\mathcal{L}(H)_+$, the notation $A\ge 0$ is also used to indicate that $A$ is a positive operator on $H$.  By a projection, we  mean an idempotent and self-adjoint element of certain $C^*$-algebra.

Throughout the rest of this paper, $\mathbb{C}^{m\times n}$ is the set of $m\times n$ complex matrices, $I_n$ is the identity matrix in $\mathbb{C}^{n\times n}$, $\mathfrak{A}$ is a $C^*$-algebra, $E$, $H$ and $K$ are Hilbert $\mathfrak{A}$-modules. The notations of ``$\oplus$" and ``$\dotplus$" are used
with different meanings. For Hilbert $\mathfrak{A}$-modules
$H_1$ and $H_2$, let
$$H_1\oplus H_2:=\left\{(h_1, h_2)^T :h_i\in H_i, i=1,2\right\},$$ which is also a Hilbert $\mathfrak{A}$-module whose $\mathfrak{A}$-valued inner product is given by
$$\left<(x_1, y_1)^T, (x_2, y_2)^T\right>=\big<x_1,x_2\big>+\big<y_1,
y_2\big>$$ for $x_i\in H_1, y_i\in H_2, i=1,2$. If both
$H_1$ and $H_2$ are submodules of a Hilbert $\mathfrak{A}$-module $H$ such that $H_1\cap H_2=\{0\}$, then we
set
$$H_1\dotplus H_2:=\{h_1+h_2 : h_i\in H_i, i=1,2\}.$$

Recall that a closed submodule $M$ of  $H$ is said to be
orthogonally complemented in $H$ if $H=M\dotplus M^\bot$, where
$$M^\bot=\big\{x\in H: \langle x,y\rangle=0\quad \mbox{for every}\ y\in
M\big\}.$$
In this case, the projection from $H$ onto $M$ is denoted by $P_M$.

An element $U$ of $\mathcal{L}(H,K)$ is said to be a partial isometry if $U^*U$ is a projection on $H$, or equivalently, $UU^*U=U$  \cite[Lemma~2.1]{XF}.
When  $H=K$, let $[A,B]=AB-BA$ be the commutator of $A,B\in\mathcal{L}(H)$.

Given $T\in\mathcal{L}(H,K)$, the polar decomposition of $T$ can be represented as
\begin{equation}\label{equ:two conditions of polar decomposition}T=U|T|\quad \mbox{and}\quad U^*U=P_{\overline{\mathcal{R}(T^*)}},\end{equation}
where $U\in\mathcal{L}(H,K)$ is a partial isometry. From \cite[Lemma~3.6 and Theorem~3.8]{LLX-AIOT}, it is known that
$T$ has the polar decomposition represented by \eqref{equ:two conditions of polar decomposition} if and only if
$\overline{\mathcal{R}(T^*)}$ and $\overline{\mathcal{R}(T)}$ are orthogonally complemented in $H$ and $K$, respectively.
In such case, the polar decomposition of $T^*$ exists and can be represented by
\begin{align}\label{equ:the polar decomposition of T star-pre stage}&T^*=U^*|T^*| \quad \mbox{and}\quad  UU^*=P_{\overline{\mathcal{R}(T)}}.
\end{align}
It is known that for the given $T\in\mathcal{L}(H,K)$, there exists at most a partial isometry $U\in\mathcal{L}(H,K)$  satisfying \eqref{equ:two conditions of polar decomposition} \cite[Lemma~3.9]{LLX-AIOT}, and such a partial isometry $U$ may however fail to be existent \cite[Example~3.15]{LLX-AIOT}. In the case that $T$ has the polar decomposition represented by \eqref{equ:two conditions of polar decomposition},  in what follows we simply say that $T$ has the polar decomposition $T=U|T|$.

The following lemmas are useful, which will be used in the sequel.

\begin{lemma}\label{lem:rang characterization-1}{\rm \cite[Proposition~2.7]{LLX-AIOT}} Let $A\in\mathcal{L}(H,K)$ and $B,C\in\mathcal{L}(E,H)$ be such that $\overline{\mathcal{R}(B)}=\overline{\mathcal{R}(C)}$. Then $\overline{\mathcal{R}(AB)}=\overline{\mathcal{R}(AC)}$.
\end{lemma}

\begin{lemma}\label{lem:Range closure of TT and T} {\rm\cite[Proposition 3.7]{Lance}}
Let $T\in\mathcal{L}(H,K)$. Then $\overline {\mathcal{R}(T^*T)}=\overline{ \mathcal{R}(T^*)}$ and $\overline {\mathcal{R}(TT^*)}=\overline{ \mathcal{R}(T)}$.
\end{lemma}

\begin{lemma}\label{lem:Range Closure of T alpha and T}{\rm (\cite[Proposition~2.9]{LLX-AIOT} and \cite[Lemma 2.2]{Vosough-Moslehian-Xu})} Let $T\in \mathcal{L}(H)_+$. Then for every $\alpha>0$, $\overline{\mathcal{R}(T^{\alpha})}=\overline{\mathcal{R}(T)}$ and $\mathcal{N}(T^{\alpha})=\mathcal{N}(T)$.
\end{lemma}

 \section{The relationship between the polar decompositions}\label{sec:equivalent characterization}
The purpose of this section is to give a full generalization of \cite[Lemma~3.1]{LLX-BJMA}. For every $T\in\mathcal{L}(H,K)$, by Lemmas~\ref{lem:Range Closure of T alpha and T} and \ref{lem:Range closure of TT and T} we have
\begin{equation}\label{equ:useful closures of range-equality} \overline{\mathcal{R}(|T|)}=\overline{\mathcal{R}(T^*T)}=\overline{\mathcal{R}(T^*)},\quad \overline{\mathcal{R}(|T^*|)}=\overline{\mathcal{R}(TT^*)}=\overline{\mathcal{R}(T)}.
\end{equation}

Our first result in this section reads as follows.

\begin{lemma}\label{lem:prepare for polar decomposition of product of operator-1}
 Suppose that $T,S\in\mathcal{L}(H)$ have the polar decompositions $T=U_T|T|$ and  $S=U_S|S|$. For every $A\in\mathcal{L}(H)$, let $X$ and $Y$ be defined by \eqref{equ:defn of X and Y}. Then
\begin{align}\label{equ:prepare for polar decomposition of product of operator-1}&|X|=U_S^*\cdot |Y|\cdot U_S,\\
\label{equ:prepare for polar decomposition of product of operator-2}&|Y|=U_S\cdot |X|\cdot U_S^*.
\end{align}
\end{lemma}
\begin{proof}Since $|Y|\in\mathcal{L}(H)_+$, by  \eqref{equ:useful closures of range-equality} we have
\begin{align}\label{equ:range Y included}\overline{\mathcal{R}(|Y|)}=\overline{\mathcal{R}(Y^*)}\subseteq \overline{\mathcal{R}(|S^*|)}=\overline{\mathcal{R}(S)}=\mathcal{R}(U_S),
\end{align}
which gives $U_SU_S^*|Y|=|Y|$, hence $|Y|=|Y|U_SU_S^*$ by taking $*$-operation. It follows that
\begin{align*}|X|^2=&S^*A^*T^*TAS=U_S^*|S^*|\cdot A^*\cdot |T|^2\cdot A\cdot |S^*|U_S\\
=&U_S^*|Y|^2U_S=(U_S^*|Y|U_S)^2,
\end{align*}
which yields $|X|=U_S^*|Y|U_S$, since both $|X|$ and $U_S^*|Y|U_S$ are positive. Consequently,
$U_S|X|U_S^*=U_SU_S^*|Y|U_SU_S^*=|Y|$.
\end{proof}

Now, we provide the technical result of this section as follows.
\begin{lemma}\label{lem:polar decomposition of product of operator-1}
Suppose that $T,S\in\mathcal{L}(H)$ have the polar decompositions $T=U_T|T|$ and  $S=U_S|S|$. For every $A\in\mathcal{L}(H)$, let $X$ and $Y$ be defined by \eqref{equ:defn of X and Y}. Then $X$ has the polar decomposition if and only if $Y$ has the polar decomposition.
\end{lemma}
\begin{proof}  As is shown in the proof of Lemma~\ref{lem:prepare for polar decomposition of product of operator-1},
we have
\begin{equation}\label{equ:simplified-01}U_S^*U_S|X|=|X|=|X|U_S^*U_S.
\end{equation}

``$\Longleftarrow$". Suppose that $Y$ has the polar decomposition $Y=U_Y|Y|$. Let
\begin{equation}\label{equ: defn of U sub X}U_X=U_TU_YU_S.\end{equation}
According to \eqref{equ:defn of X and Y}, \eqref{equ:prepare for polar decomposition of product of operator-2}, \eqref{equ:simplified-01} and \eqref{equ: defn of U sub X}, we obtain
\begin{align*}X=&U_T |T|\cdot A\cdot |S^*| U_S=U_TYU_S=U_TU_Y|Y|U_S\\
=&U_TU_YU_S|X|U_S^*U_S=U_X|X|.
\end{align*}
Meanwhile, by the definition of $Y$ and \eqref{equ:useful closures of range-equality},  we have
$$\overline{\mathcal{R}(Y)}\subseteq \overline{\mathcal{R}(T^*)}\quad\mbox{and}\quad
\overline{\mathcal{R}(Y^*)}\subseteq \overline{\mathcal{R}(S)},$$  which yield
$$U_T^*U_TU_Y=U_Y\quad \mbox{and}\quad U_SU_S^*U_Y^*=U_Y^*.$$
It follows from  \eqref{equ: defn of U sub X} that
\begin{align}\label{pre for projection-01}&U_X^*U_X=U_S^*U_Y^*(U_T^*U_TU_Y)U_S=U_S^*U_Y^*U_YU_S,\\
&(U_X^*U_X)^2=U_S^*U_Y^*U_Y (U_S U_S^*U_Y^*)U_YU_S=U_S^*U_Y^*U_YU_S=U_X^*U_X,\nonumber
\end{align}
which shows that $U_X^*U_X$ is a projection. Utilizing \eqref{equ:useful closures of range-equality}, \eqref{equ:prepare for polar decomposition of product of operator-1}, Lemmas~\ref{lem:Range closure of TT and T} and \ref{lem:rang characterization-1},  together with \eqref{pre for projection-01},  we have
\begin{align*}\overline{\mathcal{R}(X^*)}=&\overline{\mathcal{R}(|X|)}=\overline{\mathcal{R}(U_S^*|Y|U_S)}
=\overline{\mathcal{R}\left[\left(U_S^*|Y|^\frac12\right)\left(U_S^*|Y|^\frac12\right)^*\right]}
=\overline{\mathcal{R}\left(U_S^*|Y|^\frac12\right)}\\=&\overline{\mathcal{R}(U_S^*|Y|)}=\overline{\mathcal{R}(U_S^*U_Y^*)}
=\overline{\mathcal{R}(U_S^*U_Y^*U_YU_S)}=\mathcal{R}(U_X^*U_X).
\end{align*}
Hence, \eqref{equ:two conditions of polar decomposition} is satisfied such that $U,T$ therein are replaced by $U_X$ and $X$, respectively.

``$\Longrightarrow$". Suppose that $X=U_X|X|$ is the polar decomposition. Let
\begin{equation}\label{equ: defn of U sub Y}U_Y=U_T^*U_XU_S^*.\end{equation}
According to \eqref{equ:defn of X and Y}, \eqref{equ:prepare for polar decomposition of product of operator-1}, \eqref{equ: defn of U sub Y} and \eqref{equ:range Y included}, we obtain
\begin{align*}Y=&U_T^*U_T|T|\cdot A\cdot |S^*|U_SU_S^*=U_T^*T\cdot A\cdot SU_S^*=U_T^*XU_S^*=U_T^*U_X|X|U_S^*\\
=&U_T^*U_XU_S^*|Y|U_SU_S^*=U_Y|Y|.
\end{align*}
Meanwhile, by the definition of $X$   we have $\mathcal{R}(X)\subseteq \mathcal{R}(T)$ and
$\mathcal{R}(X^*)\subseteq \mathcal{R}(S^*)$, which yield
$$U_TU_T^*U_X=U_X,\quad U_S^*U_SU_X^*=U_X^*.$$
It follows from  \eqref{equ: defn of U sub Y} that
\begin{align}\label{pre for projection-03}&U_Y^*U_Y=U_SU_X^*(U_TU_T^*U_X)U_S^*=U_SU_X^*U_XU_S^*,\\
&(U_Y^*U_Y)^2=U_SU_X^*U_X (U_S^* U_SU_X^*)U_XU_S^*=U_SU_X^*U_XU_S^*=U_Y^*U_Y,\nonumber
\end{align}
hence $U_Y^*U_Y$ is a projection. Moreover, by \eqref{equ:useful closures of range-equality}, \eqref{equ:prepare for polar decomposition of product of operator-2}, Lemmas~\ref{lem:Range closure of TT and T} and \ref{lem:rang characterization-1},  together with \eqref{pre for projection-03}, we have
\begin{align*}\overline{\mathcal{R}(Y^*)}=&\overline{\mathcal{R}(|Y|)}
=\overline{\mathcal{R}(U_S|X|U_S^*)}=\overline{\mathcal{R}\left(U_S|X|^\frac12\right)}
=\overline{\mathcal{R}(U_SU_X^*)}\\=&\overline{\mathcal{R}(U_SU_X^*U_XU_S^*)}=\mathcal{R}(U_Y^*U_Y).
\end{align*}
Hence, \eqref{equ:two conditions of polar decomposition} is satisfied such that $U,T$ therein are replaced by $U_Y$ and $Y$, respectively.
\end{proof}

\begin{lemma}\label{polar decomposition wrt rho A} For every $A\in\mathcal{L}(H,K)$, let $\rho(A)\in \mathcal{L}(K\oplus H)_{\mbox{sa}}$ be defined by
\begin{equation}\label{equ:defn of rho A}\rho(A)=\left(
                          \begin{array}{cc}
                            0 & A \\
                            A^* & 0 \\
                          \end{array}
                        \right).
                        \end{equation}
Then  the following statements are valid:
\begin{enumerate}
\item[{\rm (i)}] $A$ has the polar decomposition if and only if $\rho(A)$ has the polar decomposition;
\item[{\rm (ii)}] If $A=U_A|A|$ is the polar decomposition,  then $\rho(A)=\rho(U_A)|\rho(A)|$
is the polar decomposition.
\end{enumerate}
\end{lemma}
\begin{proof}(i) Clearly, $\rho(A)^2=\mbox{diag}(AA^*,A^*A)$, so by Lemma~\ref{lem:Range closure of TT and T} we have
\begin{align*}\overline{\mathcal{R}[\rho(A)]}=\overline{\mathcal{R}[\rho(A)^2]}=\overline{\mathcal{R}(A)}\oplus \overline{\mathcal{R}(A^*)}.
\end{align*}
Therefore, $\overline{\mathcal{R}[\rho(A)]}$ is orthogonally complemented in $K\oplus H$ if and only if $\overline{\mathcal{R}(A)}$ and  $\overline{\mathcal{R}(A^*)}$ are orthogonally complemented in $K$ and $H$, respectively. Hence, the desired equivalence follows from \cite[Lemma~3.6 and Theorem~3.8]{LLX-AIOT}.

(ii) If $A=U_A|A|$ is the polar decomposition, then a simple calculation shows that $|\rho(A)|=\mbox{diag}(|A^*|,|A|)$, and $\rho(A)=\rho(U_A)|\rho(A)|$
is the polar decomposition.
\end{proof}

Combining Lemmas~\ref{lem:polar decomposition of product of operator-1} and \ref{polar decomposition wrt rho A}, a generalization can be obtained as follows.

\begin{theorem}\label{thm:main result of 3 operators}Suppose that $T\in\mathcal{L}(K)$ and $S\in\mathcal{L}(H)$
have the polar decompositions $T=U_T|T|$ and  $S=U_S|S|$.   For every $A\in\mathcal{L}(H,K)$, let $X$ and $Y$ be defined by \eqref{equ:defn of X and Y}. Then the following statements are valid:
\begin{enumerate}
\item[{\rm (i)}] If $X=U_X|X|$ is the polar decomposition, then $Y=U_T^*U_XU_S^*|Y|$ is the polar decomposition;
\item[{\rm (ii)}] If $Y=U_Y|Y|$ is the polar decomposition, then $X=U_TU_YU_S|X|$ is the polar decomposition.
\end{enumerate}
\end{theorem}
\begin{proof}Let $\widetilde{T},\widetilde{A},\widetilde{S}\in\mathcal{L}(K\oplus H)$ be defined by
$$\widetilde{T}=\left(
                  \begin{array}{cc}
                    T & 0 \\
                    0 & S^* \\
                  \end{array}
                \right),\quad \widetilde{A}=\rho(A),\quad \widetilde{S}=\widetilde{T}^*.$$
Then
$$\widetilde{X}:=\widetilde{T}\widetilde{A}\widetilde{S}=\rho(X),\quad \widetilde{Y}:=|\widetilde{T}|\cdot \widetilde{A}\cdot |\widetilde{S}^*|=\rho(Y).$$
It is clear that the polar decompositions of $\widetilde{T}$  and $\widetilde{S}$
are given by
$$\widetilde{T}=\mbox{diag}(U_T,U_S^*)\mbox{diag}(|T|,|S^*|),\quad \widetilde{S}=\mbox{diag}(U_T^*,U_S)\mbox{diag}(|T^*|,|S|).$$

Suppose that $X=U_X|X|$ is the polar decomposition. By Lemma~\ref{polar decomposition wrt rho A}
$\widetilde{X}=\rho(U_X)|\widetilde{X}|$ is the polar decomposition, which leads by  \eqref{equ: defn of U sub Y} to the polar decomposition $\widetilde{Y}=U_{\widetilde{Y}}|\widetilde{Y}|$,
where
\begin{align*}&U_{\widetilde{Y}}=U_{\widetilde{T}}^*U_{\widetilde{X}}U_{\widetilde{S}}^*=\mbox{diag}(U_T^*,U_S)\rho(U_X)\mbox{diag}(U_T,U_S^*)=\rho(U_T^*U_XU_S^*).
\end{align*}
Consequently,  by Lemma~\ref{polar decomposition wrt rho A}  $Y$ has the polar decomposition $Y=U_Y|Y|$ such that $U_{\widetilde{Y}}=\rho(U_Y)$.  Hence, $\rho(U_Y)=\rho(U_T^*U_XU_S^*)$; or equivalently, $U_Y=U_T^*U_XU_S^*$. This completes the proof of (i).  The proof of (ii) is similar.
\end{proof}

\begin{remark}
The special case of the preceding theorem  can be found in \cite[Theorem~2.1]{Ito-IEOT} and \cite[Lemma~3.1]{LLX-BJMA}, where $H=K$ and $A\in\mathcal{L}(H)$ is the identity operator.
\end{remark}

\section{The polar decomposition of a  multiplicative perturbation}\label{sec:multiplicative perturbation}

The purpose of this section is to set up the theory for the polar decomposition of a multiplicative perturbation of an operator. Suppose that $T\in\mathcal{L}(H,K)$ has the polar decomposition $T=U_T|T|$. Let $P\in\mathcal{L}(H)$ and $Q\in\mathcal{L}(K)$ be the associated projections defined by
\begin{equation}\label{defn of P and Q}P=P_{\overline{\mathcal{R}(T^*)}}\quad\mbox{and}\quad  Q=P_{\overline{\mathcal{R}(T)}}.\end{equation}
That is, $P=U_T^*U_T$ and $Q=U_TU_T^*$. Put
\begin{equation}\label{defn of H i K i}H_1=\mathcal{R}(P),\quad H_2=\mathcal{N}(P),\quad K_1=\mathcal{R}(Q),\quad K_2=\mathcal{N}(Q). \end{equation}
Assume that $H_i\ne \{0\}$ and $K_i\ne \{0\}$ for $i=1,2$. Let $U_Q: K\to  K_1\oplus K_2$ be the unitary operator defined by
\begin{equation}\label{equ:defn of W T} U_Q x=\big(Q x, (I_K-Q)x\big)^T,\quad \forall x\in K.
\end{equation}
It is easy to verify that $U_Q^*\in \mathcal{L}(K_1\oplus K_2, K)$ is given  by
 \begin{eqnarray*} U_Q^*\big((x_1, x_2)^T\big)=x_1+x_2,\quad \forall\,x_i\in K_i, i=1,2.\end{eqnarray*}
Similarly, define the unitary operator $U_P: H\to H_1\oplus H_2$. Thus, for every $E\in\mathcal{L}(K)$ we have
\begin{equation}\label{equ:partitioned form of E}U_Q E U_Q^*=\left(
                                      \begin{array}{cc}
                                        B & E_{12} \\
                                        C & E_{22} \\
                                      \end{array}
                                    \right)\in\mathcal{L}(K_1\oplus K_2),
\end{equation}
where $B=Q E Q|_{K_1}$ is the restriction of $Q E Q$ on $K_1$. Similarly,
$C=(I_K-Q)EQ|_{K_1}$, $E_{12}=Q E(I_K-Q)|_{K_2}$ and $E_{22}=(I_K-Q)E(I_K-Q)|_{K_2}$. Specifically,
\begin{equation}\label{equ:partitioned form of Q}U_Q Q U_Q^*=\left(
                              \begin{array}{cc}
                                I_{K_1} & 0 \\
                                0 & 0 \\
                              \end{array}
                            \right).\end{equation}
Utilizing $U_Q(EQ)U_Q^*=U_QEU_Q^*\cdot U_Q Q U_Q^*$, by \eqref{equ:partitioned form of E} and \eqref{equ:partitioned form of Q} we obtain
\begin{equation}\label{equ:partitioned form of EQ}U_Q(EQ)U_Q^*=\left(
                                                          \begin{array}{cc}
                                                            B & 0 \\
                                                            C & 0 \\
                                                          \end{array}
                                                        \right).
\end{equation}
Similarly, for every $F\in\mathcal{L}(H)$, there exist $D\in\mathcal{L}(H_1)$ and $G\in\mathcal{L}(H_1,H_2)$ such that
\begin{equation}\label{equ:partitioned form of FP}U_P(FP)U_P^*=\left(
                               \begin{array}{cc}
                                 D & 0 \\
                                 G & 0 \\
                               \end{array}
                             \right).
\end{equation}
It follows from \eqref{defn of P and Q} that $QT=T=TP$, which is combined with
$$U_Q TU_P^*=\left(
                             \begin{array}{cc}
                               QTP|_{H_1} & QT(I_H-P)|_{H_2} \\
                               (I_K-Q)TP|_{H_1} & (I_K-Q)T(I_H-P)|_{H_2} \\
                             \end{array}
                           \right)$$
to get
\begin{align}\label{equ:partitioned form of T}
U_Q TU_P^*=&\left(
    \begin{array}{cc}
      A & 0 \\
      0 & 0 \\
    \end{array}
  \right)=:\widetilde{A},
\end{align}
where $A=QTP|_{H_1}\in\mathcal{L}(H_1,K_1)$. So, when $M$ is a multiplicative of $T$ given by $M=ETF^*$, we have
$M=(EQ)T(FP)^*$, and thus by \eqref{equ:partitioned form of EQ}--\eqref{equ:partitioned form of T}
\begin{align}\label{equ:partitioned form of M}\widetilde{X}:=U_Q MU_P^*=\left(
                        \begin{array}{cc}
                          B & 0 \\
                          C & 0 \\
                        \end{array}
                      \right)\left(
                               \begin{array}{cc}
                                 A & 0 \\
                                 0 & 0 \\
                               \end{array}
                             \right)\left(
                                      \begin{array}{cc}
                                        D & 0 \\
                                        G & 0 \\
                                      \end{array}
                                    \right)^*=:\widetilde{T}\widetilde{A}\widetilde{S}.
\end{align}
It is clear that
\begin{align*}|\widetilde{T}|=\mbox{diag}\big(\theta_{B,C},0\big)\quad\mbox{and}\quad |\widetilde{S}^*|=\mbox{diag}\big(\theta_{D,G},0\big),
\end{align*}
where $\theta_{B,C}\in \mathcal{L}(K_1)$ and $\theta_{D,G}\in \mathcal{L}(H_1)$ are given by
\begin{equation}\label{equ:defn of theta++}\theta_{B,C}=(B^*B+C^*C)^\frac12\quad\mbox{and}\quad  \theta_{D,G}=(D^*D+G^*G)^\frac12.\end{equation}
Therefore,
\begin{equation}\label{defn of widetilde Y}\widetilde{Y}:=|\widetilde{T}|\cdot \widetilde{A}\cdot |\widetilde{S}^*|=\mbox{diag}(Y,0),\end{equation}
in which
\begin{equation}\label{equ:concrete expression of Y}Y=\theta_{B,C}\cdot A\cdot \theta_{D,G}\in\mathcal{L}(H_1,K_1).\end{equation}

To deal with the polar decompositions of $\widetilde{T}$ and $\widetilde{S}$, we need some knowledge obtained recently in \cite{Fu etal}.

\begin{definition}\cite{Fang-Moslehian-Xu}  Let $A\in\mathcal{L}(H,K)$ and $C\in\mathcal{L}(E,K)$. An operator $D\in\mathcal{L}(E,H)$ is said to be the reduced solution
to the system
\begin{equation}\label{reduced solution system}AX=C,\quad  X\in\mathcal{L}(E,H),\end{equation} if $AD=C$ and $\mathcal{R}(D)\subseteq \overline{\mathcal{R}(A^*)}$.
\end{definition}

Since $\overline{\mathcal{R}(A^*)}\subseteq \mathcal{N}(A)^\bot$,  such a reduced solution (if it exists) is unique.

\begin{lemma}\label{lem:two orthogonally complemented conditions} {\rm \cite[Lemma~1.3]{MMX}} Let $A\in\mathcal{L}(H,K)$. Then the following statements are equivalent:
\begin{enumerate}
\item[{\rm (i)}] $\overline{{\mathcal R}(A^*)}$ is orthogonally complemented in $H$;
\item[{\rm (ii)}] For every $C\in\mathcal{L}(E, K)$ with ${\mathcal R}(C)\subseteq {\mathcal R}(A)$, system \eqref{reduced solution system}
has the reduced solution.
\end{enumerate}
\end{lemma}

\begin{definition}\label{defn of tractable}\cite[Definition~3.1]{Fu etal}
For every  $A\in\mathcal{L}(H,K)$ and $B\in\mathcal{L}(K)$, let $T_{A,B}\in \mathcal{L}(H\oplus K)$ and $S_{A,B}\in\mathcal{L}(K)$ be defined by
\begin{equation}\label{equ.1:associated operator T}
T_{A,B}=
\left(\begin{array}{cc}
0 & 0 \\
A & B \\
\end{array}\right)\quad\mbox{and}\quad  S_{A,B}=(AA^*+BB^*)^\frac12.
\end{equation}
\end{definition}

\begin{lemma}\label{lem:characterization of the tractable pair}{\rm \cite[Theorem~3.1]{Fu etal}} For every $A\in\mathcal{L}(H,K)$ and $B\in\mathcal{L}(K)$, let  $T_{A,B}$ and $S_{A,B}$ be defined by \eqref{equ.1:associated operator T}. Then $T_{A,B}$ has the polar decomposition
if and only if the following conditions are satisfied:
\begin{itemize}
\item[{\rm (i)}] $\overline{\mathcal{R}(S_{_{A,B}})}$ is orthogonally complemented in $K$;
\item[{\rm (ii)}] $\mathcal{R}(A)\subseteq \mathcal{R}(S_{_{A,B}})$ and  $\mathcal{R}(B)\subseteq \mathcal{R}(S_{_{A,B}})$.
\end{itemize}
\end{lemma}

\begin{remark}\label{rem:polar decomposition for 2 by 2 block matrices} Suppose that conditions (i) and (ii) in Lemma~\ref{lem:characterization of the tractable pair} are satisfied. It follows from the proof of \cite[Theorem~3.1]{Fu etal} that $T^*_{A,B}=U^*|T^*_{A,B}|$ is the polar decomposition,  where
$$U=\left(
  \begin{array}{cc}
    0 & 0 \\
    C & D \\
  \end{array}
\right)\quad\mbox{and}\quad  |T^*_{A,B}|=\left(
\begin{array}{cc}
0 & 0 \\
0 & S_{A,B} \\
\end{array}
\right),$$
in which $C$ and $D$ are the reduced solutions of
\begin{align*}
S_{A,B}X=A\ \big(X\in\mathcal{L}(H,K)\big)\quad\mbox{and}\quad S_{A,B}Y=B \ \big(Y\in\mathcal{L}(K)\big),
\end{align*}
respectively.
\end{remark}

Now, we provide the main result of this section as follows.
\begin{theorem}\label{thm:polar decomposition for multiplicative perturbation} Suppose that $T\in\mathcal{L}(H,K)$ has  the polar decomposition such that $H_i, K_i (i=1,2)$ defined by \eqref{defn of H i K i} are nonzero. Let $M=ETF^*$ be a multiplicative perturbation of $T$ with $E\in \mathcal{L}(K)$ and $F\in \mathcal{L}(H)$,
and let $P,Q$, $B,C$, $D,G$, $A$, $\widetilde{A}$, $\widetilde{X}$, $\theta_{B,C}$, $\theta_{D,G}$, $\widetilde{Y}$ and $Y$  be defined by \eqref{defn of P and Q} and \eqref{equ:partitioned form of EQ}--\eqref{equ:concrete expression of Y}, respectively.
Suppose that the following conditions are  satisfied:
\begin{enumerate}
\item[{\rm (i)}] $\overline{\mathcal{R}(\theta_{B,C})}$ and $\overline{\mathcal{R}(\theta_{D,G})}$
 are orthogonally complemented in $K_1$ and $H_1$, respectively.
\item[{\rm (ii)}] $\mathcal{R}(B^*)\subseteq \mathcal{R}(\theta_{B,C}), \mathcal{R}(C^*)\subseteq \mathcal{R}(\theta_{B,C})$,  $\mathcal{R}(D^*)\subseteq \mathcal{R}(\theta_{D,G})$ and  $\mathcal{R}(G^*)\subseteq \mathcal{R}(\theta_{D,G})$.
\item[{\rm (iii)}] $\overline{\mathcal{R}(Y^*)}$ and $\overline{\mathcal{R}(Y)}$ are orthogonally complemented in $H_1$ and $K_1$, respectively.
\end{enumerate}
Then   $\widetilde{X}=U_{\widetilde{X}}|\widetilde{X}|$ is the polar decomposition such that
 \begin{align}\label{formula for U widetilde X} U_{\widetilde{X}}=\left(
                                                                    \begin{array}{cc}
                                                                      Z_1U_YZ_3^* & Z_1U_YZ_4^* \\
                                                                      Z_2U_YZ_3^* & Z_2U_YZ_4^* \\
                                                                    \end{array}
                                                                  \right),
\end{align}
where $Y=U_Y|Y|$ is the polar decomposition,
 $Z_1^*$, $Z_2^*$, $Z_3^*$ and $Z_4^*$ are the reduced solutions of
\begin{align}\label{reduced system-01}&\theta_{B,C}X_1=B^*\ \big(X_1\in\mathcal{L}(K_1)\big),\quad \theta_{B,C}X_2=C^* \ \big(X_2\in\mathcal{L}(K_2,K_1)\big),\\
\label{reduced system-02}&\theta_{D,G}X_3=D^*\ \big(X_3\in\mathcal{L}(H_1)\big),\quad \theta_{D,G}X_4=G^* \ \big(X_4\in\mathcal{L}(H_2,H_1)\big),
\end{align}
respectively.
\end{theorem}
\begin{proof}By \eqref{equ:partitioned form of M},   Lemma~\ref{lem:characterization of the tractable pair} and Remark~\ref{rem:polar decomposition for 2 by 2 block matrices}, the polar decompositions of $\widetilde{T}$ and $\widetilde{S}^*$ are given by $\widetilde{T}=U_{\widetilde{T}}|\widetilde{T}|$ and  $\widetilde{S}^*=U_{\widetilde{S}^*}|\widetilde{S}^*|$ respectively, where
\begin{align*}&U_{\widetilde{T}}=\left(
                                   \begin{array}{cc}
                                     Z_1 & 0 \\
                                     Z_2 & 0 \\
                                   \end{array}
                                 \right),\quad |\widetilde{T}|=\left(
                                                                 \begin{array}{cc}
                                                                   \theta_{B,C} & 0 \\
                                                                   0 & 0 \\
                                                                 \end{array}
                                                               \right),\\
&U_{\widetilde{S}^*}=\left(
                                   \begin{array}{cc}
                                     Z_3 & 0 \\
                                     Z_4 & 0 \\
                                   \end{array}
                                 \right),\quad |\widetilde{S}^*|=\left(
                                                                 \begin{array}{cc}
                                                                   \theta_{D,G} & 0 \\
                                                                   0 & 0 \\
                                                                 \end{array}
                                                               \right).
\end{align*}
By assumption $Y$  has the polar decomposition $Y=U_Y|Y|$, so from \eqref{defn of widetilde Y} we can obtain
$ U_{\widetilde{Y}}=\mbox{diag}(U_Y,0)$.
Thus, according to Theorem~\ref{thm:main result of 3 operators}~(ii) the polar decomposition of $\widetilde{X}$ is given by
 $\widetilde{X}=U_{\widetilde{X}}|\widetilde{X}|$, where
\begin{align*}&U_{\widetilde{X}}=U_{\widetilde{T}}U_{\widetilde{Y}}U_{\widetilde{S}}=\left(
                                                                                       \begin{array}{cc}
                                                                                         Z_1 & 0 \\
                                                                                         Z_2 & 0 \\
                                                                                       \end{array}
                                                                                     \right)\left(
                                                                                              \begin{array}{cc}
                                                                                                U_Y & 0 \\
                                                                                                0 & 0 \\
                                                                                              \end{array}
                                                                                            \right)\left(
                                                                                                     \begin{array}{cc}
                                                                                                       Z_3^* & Z_4^* \\
                                                                                                       0 & 0 \\
                                                                                                     \end{array}
                                                                                                   \right),
\end{align*}
which can be simplified to \eqref{formula for U widetilde X}, as desired.
\end{proof}

\begin{remark}It is remarkable that conditions (i)--(iii) in Theorem~\ref{thm:polar decomposition for multiplicative perturbation} are satisfied automatically in the Hilbert space case.
\end{remark}

\begin{remark}The special case of Theorem~\ref{thm:polar decomposition for multiplicative perturbation} can be found in the proof of \cite[Lemma~2.3]{Fu-Xu}, where  $M=ETE^*$, in which $T$ is a Hermitian matrix, and $E$ ia a square matrix such that $\theta_{B,C}$ defined by \eqref{equ:defn of theta++} is nonsingular.
\end{remark}

\section{The polar decomposition of the product of three operators}\label{sec:three operators}
The purpose of this section is to give new characterizations of the polar decomposition associated with three operators. Recall that an operator $W$ is said to be a contraction if $\|W\|\le 1$. A useful lemma concerning contractions can be provided as follows.

\begin{lemma}\label{lem:trivial add useful}Let $T,W\in\mathcal{L}(H)$ be given such that $W$ is a contraction.
 \begin{enumerate}
 \item[{\rm (i)}] If $T=W|T|$, then $|T|=W^*T$. Furthermore, $T\ge 0$ whenever $W$ is a projection.
 \item[{\rm (ii)}] If $|T|=W T$, then $T=W^*|T|$. Furthermore, $T\ge 0$ whenever $W$ is a projection.
 \end{enumerate}
\end{lemma}
\begin{proof}(i) By assumption we have
$|T|^2=|T|W^*W|T|$, hence
$$|T|(I_H-W^*W)|T|=0.$$
Since $\|W\|\le 1$, the equation above can be rephrased as
$$\left[(I_H-W^*W)^\frac12 |T|\right]^*
(I_H-W^*W)^\frac12 |T|=0,$$
 which implies that
 \begin{equation}\label{equ:noname-01}|T|=W^*W|T|,\end{equation}
and thus $|T|=W^*T$, since $T=W|T|$. If furthermore $W$ is a projection, then $W^*W=W$, so by \eqref{equ:noname-01} $|T|=W|T|=T$, hence $T\ge 0$.

 (ii) By assumption $T^*(I_H-W^*W)T=0,$
 which yields $T=W^*WT=W^*|T|$. If furthermore $W$ is a projection, then $T=W^*WT=WT=|T|$, hence $T\ge 0$.
\end{proof}

We provide an auxiliary lemma as follows.
\begin{lemma}Let $T,A,S\in\mathcal{L}(H)$ be such that $T, S$ and $X$ defined by \eqref{equ:defn of X and Y}  have the polar decompositions
$T=U_T|T|, S=U_S|S|$ and $X=U_X|X|$. Then there exists a contraction $B\in\mathcal{L}(H)$ such that $U_X=U_TBU_S$.
\end{lemma}
\begin{proof}Since $\mathcal{R}(U_X)=\overline{\mathcal{R}(X)}\subseteq \overline{\mathcal{R}(T)}$ and $\mathcal{R}(U_X^*)=\overline{\mathcal{R}(X^*)}\subseteq \overline{\mathcal{R}(S^*)}$, we have $U_TU_T^*U_X=U_X$ and $U_S^*U_S U_X^*=U_X^*$. Hence, $U_X=U_TBU_S$ with $B=U_T^*U_XU_S^*$.
\end{proof}
We provide a technical result as follows.
\begin{lemma}\label{lem:technical lemma for three operators} Suppose that $T,A,S,B\in\mathcal{L}(H)$ are such that $T$ and $S$ have the polar decompositions
$T=U_T|T|$ and $S=U_S|S|$, and $B$ is a contraction. Let $X$ and $Y$ be defined by \eqref{equ:defn of X and Y} and let $W=U_T^*U_T BU_SU_S^*$.
Then the  following statements are  equivalent:
\begin{enumerate}
\item[{\rm (i)}]$X=U_TBU_S|X|$;
\item[{\rm (ii)}]$Y=W|Y|$;
\item[{\rm (iii)}] $B^*Y\ge 0$ and $BB^*Y=Y$.
\end{enumerate}
\end{lemma}
\begin{proof}By definition we have
\begin{equation}\label{defn of Y and W}X=TAS,\quad Y=|T|\cdot A\cdot |S^*|\quad \mbox{and}\quad W=W_1U_SU_S^*,\end{equation}
in which $W_1=U_T^*U_TB$. It is clear that
\begin{equation}\label{equ:equation wrt Y and abs of Y}U_T^*U_T Y=YU_SU_S^*=Y, \quad U_SU_S^*|Y|=|Y|, \quad U_T^*U_TW=W.\end{equation}

(i)$\Longrightarrow$(ii). By \eqref{equ:prepare for polar decomposition of product of operator-1} we have
\begin{align*}U_T Y U_S=X=U_TBU_S U_S^*|Y|U_S,
\end{align*}
which is combined with the first equation in \eqref{equ:equation wrt Y and abs of Y} to get
\begin{align*}Y=U_T^*(U_TYU_S)U_S^*=W|Y|U_SU_S^*=W|Y|.
\end{align*}

(ii)$\Longrightarrow$(i). It follows from \eqref{equ:prepare for polar decomposition of product of operator-2} that
\begin{align*}X=U_T Y U_S=U_TW U_S|X|U_S^*U_S=U_TW U_S|X|=U_TBU_S|X|.
\end{align*}

(ii)$\Longrightarrow$(iii). Suppose that $Y=W|Y|$. Then $Y=W_1U_SU_S^*|Y|=W_1|Y|$, hence by Lemma~\ref{lem:trivial add useful}
$W_1^*Y=|Y|$, which yields $B^* Y=|Y|\ge 0$, since $W_1^*Y=B^* Y$.
So, by Lemma~\ref{lem:trivial add useful} $Y=B|Y|$. Thus,  $BB^*Y=B|Y|=Y$.

(iii)$\Longrightarrow$(ii). Suppose that $B^* Y\ge 0$ and $BB^* Y=Y$. Then
\begin{equation}\label{another form of abs Y}(B^* Y)^2=(B^* Y)^*(B^* Y)=Y^*BB^*Y=Y^*Y=|Y|^2,\end{equation}
hence $B^* Y=|Y|$, so  $Y=B|Y|$ by  Lemma~\ref{lem:trivial add useful}, which in turn leads by \eqref{equ:equation wrt Y and abs of Y} to
$$Y=U_T^*U_TBU_SU_S^*|Y|=W|Y|.$$
This shows the equivalence of (i), (ii) and (iii).
\end{proof}

Now, we give a characterization of the polar decomposition of the product of three operators as follows.
\begin{theorem}\label{thm:technical lemma for three operators} Under the notations and conditions of Lemma~\ref{lem:technical lemma for three operators}, the  following statements are  equivalent:
\begin{enumerate}
\item[{\rm (i)}]$X=U_TBU_S|X|$ is the polar decomposition;
\item[{\rm (ii)}]$Y=W|Y|$ is the polar decomposition;
\item[{\rm (iii)}] $B^*Y\ge 0$, $BB^*Y=Y$ and $\overline{\mathcal{R}(W)}=\overline{\mathcal{R}(Y)}$;
\item[{\rm (iv)}] $B^*Y\ge 0$, $BB^*Y=Y$  and $\overline{\mathcal{R}(W^*)}=\overline{\mathcal{R}(Y^*)}$.
\end{enumerate}
\end{theorem}
\begin{proof} (i)$\Longleftrightarrow$(ii). Let $U_X=U_TBU_S$. It is clear that $W=U_T^*U_XU_S^*$ and $U_TWU_S=U_X$.
Hence, the equivalence of (i) and (ii) is immediate from Theorem~\ref{thm:main result of 3 operators}.

(ii)$\Longrightarrow$(iii) and (ii)$\Longrightarrow$(iv). From Lemma~\ref{lem:technical lemma for three operators}, we obtain $B^*Y\ge 0$ and $BB^*Y=Y$. Furthermore, by \eqref{equ:two conditions of polar decomposition} and \eqref{equ:the polar decomposition of T star-pre stage} we have $\overline{\mathcal{R}(W^*)}=\overline{\mathcal{R}(Y^*)}$ and $\overline{\mathcal{R}(W)}=\overline{\mathcal{R}(Y)}$.

(iii)$\Longrightarrow$(ii).  Since $B^*Y\ge 0$, we have
$$B^*Y=(B^*Y)^*=Y^* B=|S^*|\cdot (A^*\cdot |T|\cdot B),$$ hence
$U_SU_S^* B^*Y=B^*Y$, which means that
\begin{equation}\label{temp-03}U_SU_S^* B^*x=B^*x,\quad\forall\,x\in \overline{\mathcal{R}(Y)}.\end{equation}
By assumption $\overline{\mathcal{R}(W)}=\overline{\mathcal{R}(Y)}$, so from \eqref{temp-03} we obtain
\begin{equation}\label{temp-01}U_SU_S^* B^*W=B^*W.\end{equation}
 The same argument applied to the assumption $BB^* Y=Y$ gives $BB^*W=W$. This, together with \eqref{temp-01} and the last equation in \eqref{equ:equation wrt Y and abs of Y} yields
\begin{align*}WW^*W=&U_T^*U_TBU_SU_S^*B^*W=U_T^*U_TBB^*W=U_T^*U_TW=W.
\end{align*}
This shows that $W$ is a partial isometry. Furthermore, in virtue of Lemma~\ref{lem:rang characterization-1}, \eqref{temp-01}, the first equation in \eqref{equ:range Y included}, and \eqref{another form of abs Y},  we have
\begin{align*}\mathcal{R}(W^*)=&\mathcal{R}(W^*W)=\overline{\mathcal{R}(W^*W)}=\overline{\mathcal{R}(U_SU_S^*B^*W)}=\overline{\mathcal{R}(B^*W)}
\\=&\overline{\mathcal{R}(B^*Y)}=\overline{\mathcal{R}(|Y|)}=\overline{\mathcal{R}(Y^*)}.
\end{align*}
This completes the proof that $Y=W|Y|$ is the polar decomposition.

(iv)$\Longrightarrow$(iii). Since $B^*Y\ge 0$ and $BB^*Y=Y$, we have $Y=W|Y|$ by Lemma~\ref{lem:technical lemma for three operators}. This, together with $\overline{\mathcal{R}(W^*)}=\overline{\mathcal{R}(Y^*)}$, Lemmas~\ref{lem:rang characterization-1} and \ref{lem:Range closure of TT and T}, yields
$$\overline{\mathcal{R}(Y)}=\overline{\mathcal{R}(W|Y|)}=\overline{\mathcal{R}(WW^*)}=\overline{\mathcal{R}(W)}.\qedhere$$
\end{proof}

\begin{remark}Suppose that $T,A,S\in\mathcal{L}(H)$ are such that $T, A$ and $S$ have the polar decompositions
$T=U_T|T|$, $A=U_A|A|$ and $S=U_S|S|$. Let $X$ and $Y$ be defined by \eqref{equ:defn of X and Y} such that $X$ has the polar decomposition $X=U_X|X|$.
Naturally, one may wonder whether it is true that $U_X=U_TU_AU_S$? The following Example~\ref{ex:not natrual} shows that such an equation may be false even if conditions (i)--(iii) stated in Lemma~\ref{lem:technical lemma for three operators} are satisfied.
\end{remark}

\begin{example}\label{ex:not natrual}Let $H=\mathbb{C}^3$, and let $T,A,S\in\mathcal{L}(H)$ be given by
\begin{align*} T = \left(
                                                                    \begin{array}{ccc}
                                                                      1 & 0 & 0 \\
                                                                      0 & 1 & 0 \\
                                                                      0 & 0 & 0 \\
                                                                    \end{array}
                                                                  \right),\quad
A = \left(
                                                                    \begin{array}{ccc}
                                                                      1 & 0 & 0 \\
                                                                      0 & 0 & 0 \\
                                                                      1 & 0 & 1 \\
                                                                    \end{array}
                                                                  \right),\quad
S = \left(
                                                                    \begin{array}{ccc}
                                                                      0 & 0 & 0 \\
                                                                      0 & 1 & 0 \\
                                                                      0 & 0 & 1 \\
                                                                    \end{array}
                                                                  \right).
\end{align*}
Obviously, $U_T=U_T^*=|T|=T$ and $U_S=U_S^*=|S^*|=S$. It is routine to verify that $A=U_A|A|$ is the polar decomposition, where
\begin{align*}
U_A = \left(
                                                                    \begin{array}{ccc}
                                                                      \frac{2\sqrt{5}}{5} & 0 & -\frac{\sqrt{5}}{5} \\
                                                                      0 & 0 & 0 \\
                                                                      \frac{\sqrt{5}}{5} & 0 & \frac{2\sqrt{5}}{5} \\
                                                                    \end{array}
                                                                  \right)\quad\mbox{and}\quad
|A|=\left(
                                                                    \begin{array}{ccc}
                                                                      \frac{3\sqrt{5}}{5} & 0 & \frac{\sqrt{5}}{5} \\
                                                                      0 & 0 & 0 \\
                                                                      \frac{\sqrt{5}}{5} & 0 & \frac{2\sqrt{5}}{5} \\
                                                                    \end{array}
                                                                  \right).
\end{align*}
Let $Y$ and $W$ be defined by \eqref{defn of Y and W}. Direct computations yield
\begin{align*} Y=\left(
                                                                    \begin{array}{ccc}
                                                                      0 & 0 & 0 \\
                                                                      0 & 0 & 0 \\
                                                                      0 & 0 & 0 \\
                                                                    \end{array}
                                                                  \right)\quad\mbox{and}\quad
W=\left(
                                                                    \begin{array}{ccc}
                                                                      0 & 0 & -\frac{\sqrt{5}}{5} \\
                                                                      0 & 0 & 0 \\
                                                                      0 & 0 & 0 \\
                                                                    \end{array}
                                                                  \right).
\end{align*}
Since $Y=0$ and $W\ne 0$, it is clear that $U_A^*Y\ge 0$, $U_AU_A^*Y=Y$, $\overline{\mathcal{R}(W)}\ne \overline{\mathcal{R}(Y)}$ and $\overline{\mathcal{R}(W^*)}\ne \overline{\mathcal{R}(Y^*)}$. Thus, by Theorem~\ref{thm:polar decom of three operators} $X=U_TU_AU_S|X|$ fails to be the polar decomposition.
\end{example}

\begin{remark}Following the notations as in Lemma~\ref{lem:technical lemma for three operators}, we show that there exist Hilbert $C^*$-module $H$,  and $T,A\in\mathcal{L}(H)$ such that $X=U_TU_AU_S|X|$ with $T=S$ is the polar decomposition, and the following statements are valid:
\begin{enumerate}
\item[{\rm (i)}] $A\ne 0$ and  $AT=TA$;
\item[{\rm (ii)}] $T^{2n-1}=U_T^{2n-1}|T^{2n-1}|$ is the polar decomposition for every $n\in\mathbb{N}$;
\item[{\rm (iii)}] $T^{2n}\ne U_T^{2n}|T^{2n}|$ for every $n\in \mathbb{N}$;
\item[{\rm (iii)}] $T^{n}A=U_T^{n}U_A|T^{n}A|$ is the polar decomposition for every $n\in\mathbb{N}$.
\end{enumerate}
\end{remark}

\begin{example}Let $H=\mathbb{C}^3$, and let $T,A,S\in\mathcal{L}(H)$ be given by
\begin{align*} T = S=\left(
                                         \begin{array}{ccc}
                                           2 & 0 & 1 \\
                                           0 & 1 & 0 \\
                                           0 & 0 & -1 \\
                                         \end{array}
                                       \right),\quad
A = \left(
                                                                    \begin{array}{ccc}
                                                                      0 & 0 & 0 \\
                                                                      0 & 1 & 0 \\
                                                                      0 & 0 & 0 \\
                                                                    \end{array}
                                                                  \right).
\end{align*}
It is clear that $U_A=|A|=TA=AT=A$. Hence, $X:=TAS=A$ and $X=A|X|$ is the polar decomposition.
Direct computations yield the polar decomposition
$T=U_T|T|$, where
\begin{eqnarray}\label{equ:defn of a and U}
a=\frac{\sqrt{10}}{10}, \quad U_T=\left(
                       \begin{array}{ccc}
                        3a & 0 & a \\
                       0 & 1 & 0 \\
                      a & 0 & -3a \\
                   \end{array}
                \right),\quad   |T|=\left( \begin{array}{ccc}
                                           6a & 0 & 2a \\
                                           0 & 1 & 0 \\
                                           2a & 0 & 4a \\
                                         \end{array}
                                       \right).
\end{eqnarray}
It follows that $U_TU_AU_S=U_TAU_T=A$. This shows that $X=U_TU_AU_S|X|$ is the polar decomposition.

A simple calculation shows that $T=PDP^{-1}$, where
\begin{eqnarray*}P=\left(
                       \begin{array}{ccc}
                         1 & 0 & 1 \\
                         0 & 1 & 0 \\
                         0 & 0 & -3 \\
                       \end{array}
                     \right),\quad D=\left(
                                    \begin{array}{ccc}
                                      2 & 0 & 0 \\
                                      0 & 1 & 0 \\
                                      0 & 0 & -1 \\
                                    \end{array}
                                  \right),\quad
                                  P^{-1}=\left(
                                    \begin{array}{ccc}
                                      1 & 0 & \frac{1}{3} \\
                                      0 & 1 & 0 \\
                                      0 & 0 & -\frac{1}{3} \\
                                    \end{array}
                                  \right),
\end{eqnarray*}
which leads to
\begin{equation}\label{equ:expression of T power k-concrete}
T^{k}=PD^{k}P^{-1}=\left(
                       \begin{array}{ccc}
                         2^{k} & 0 & \frac{2^{k}+(-1)^{k+1}}{3}\\
                         0 & 1 & 0 \\
                         0 & 0 & (-1)^{k} \\
                       \end{array}
                     \right),\quad\forall k\in\mathbb{N}.
\end{equation}
Consequently, for every $n\in\mathbb{N}$
\begin{equation}\label{equ:expression of T 2n-1}
T^{2n-1}=\left(
                       \begin{array}{ccc}
                         b & 0 & c \\
                         0 & 1 & 0 \\
                         0 & 0 & -1 \\
                       \end{array}
                     \right), \quad\mbox{where $b=2^{2n-1}$, $c=\frac{1}{3}\cdot\big(2^{2n-1}+1\big)$}.
\end{equation}
Therefore,
\begin{equation*}(T^{2n-1})^*\cdot T^{2n-1}=\left(
                       \begin{array}{ccc}
                         b^2 & 0 & bc\\
                         0 & 1 & 0 \\
                         bc & 0 & 1+c^2\\
                       \end{array}
                     \right)=V\left(
                                \begin{array}{ccc}
                                  1 & 0 & 0 \\
                                  0 & \lambda_1^2 & 0 \\
                                  0 & 0 & \lambda_2^2 \\
                                \end{array}
                              \right)V^*,\end{equation*}
in which
$$p=(b+1)^2+c^2, \quad q=(b-1)^2+c^2, \quad \lambda_1=\frac{\sqrt{p}+\sqrt{q}}{2}, \quad\lambda_2=\frac{\sqrt{p}-\sqrt{q}}{2},$$
and $V$ is a unitary matrix given by
\begin{equation*}V=\left(
                 \begin{array}{ccc}
                      0 & \frac{\sqrt{b^2-\lambda_2^2}}{(pq)^\frac14} & \frac{\sqrt{\lambda_1^2-b^2}}{(pq)^\frac14} \\
                       1 & 0 & 0 \\
                      0 & \frac{\sqrt{\lambda_1^2-b^2}}{(pq)^\frac14} & -\frac{\sqrt{b^2-\lambda_2^2}}{(pq)^\frac14} \\
                    \end{array}
                 \right).
\end{equation*}
It follows that
\begin{equation*}\label{equ:half T 2n-1}\big|T^{2n-1}\big|=V\left(
                                \begin{array}{ccc}
                                  1 & 0 & 0 \\
                                  0 & \lambda_1 & 0 \\
                                  0 & 0 & \lambda_2 \\
                                \end{array}
                              \right)V^*=\left(
                      \begin{array}{ccc}
                         2^{2n-1}\cdot (3a) & 0 & 2^{2n-1}\cdot a \\
                         0 & 1 & 0 \\
                         2^{2n-1}\cdot a & 0 & \frac{2^{2n-1}+10}{3}\cdot  a \\
                      \end{array}
                    \right),
\end{equation*}
which is combined with \eqref{equ:expression of T 2n-1} and  \eqref{equ:defn of a and U} to get
$$T^{2n-1}\big|T^{2n-1}\big|^{-1}=\left(
                       \begin{array}{ccc}
                        3a & 0 & a \\
                       0 & 1 & 0 \\
                      a & 0 & -3a \\
                   \end{array}
                \right)=U_T.$$
Since $T$ is invertible, this shows that $T^{2n-1}=U_T|T^{2n-1}|$ is the polar decomposition for every $n\in\mathbb{N}$.
From the expression of $U_T$ given as above, we obtain $U_T^2=I_3$. Hence, $U_T^{2n-1}=U_T$ for every $n\in\mathbb{N}$, and thus $T^{2n-1}=U_T^{2n-1}|T^{2n-1}|$ is the polar decomposition.

For every $n\in\mathbb{N}$, by \eqref{equ:expression of T power k-concrete} it is clear that $T^{2n}\ne \big(T^{2n}\big)^*$, which means
that $T^{2n}\ne U_T^{2n}|T^{2n}|$, since $U_T^{2}=I_3$, as is shown above.  Finally, a repeated use of $TA=A$ and $U_TU_A=A$ gives $T^nA=A$ and $U_T^nU_A=A$ for every $n\in\mathbb{N}$. This shows that $T^{n}A=U_T^{n}U_A|T^{n}A|$ is the polar decomposition.

\end{example}

To deal with applications of  Theorem~\ref{thm:polar decom of three operators}, we need some knowledge about the commutators of operators.

\begin{definition}\label{def:commutator of a and b}
 For every $T,S\in\mathcal{L}(H)$, let $[T,S]=TS-ST$  be the commutator of $T$ and $S$.
\end{definition}

\begin{lemma}\label{lem:commutative property extended-0}{\rm \cite[Lemma~2.2~(ii)]{LLX-BJMA}} For every $T,S\in \mathcal{L}(H)_+$, $TS\ge 0$ if and only if $[T,S]=0$.
\end{lemma}

\begin{lemma}\label{lem:commutative property extended-3}{\rm \cite[Propositions~2.4 and 2.6]{LLX-AIOT} Let $S\in \mathcal{L}(H)$ and $T\in\mathcal{L}(H)_+$ be such that $[S,T]=0$. Then the following statements are valid:
\begin{enumerate}
\item[{\rm (i)}] $[S,T^\alpha]=0$ for every $\alpha>0$;
\item[{\rm (ii)}] If $\overline{\mathcal{R}(T)}$ is orthogonally complemented in $H$, then $\big[S,  P_{\overline{\mathcal{R}(T)}}\big]=0$.
\end{enumerate}
}\end{lemma}

\begin{lemma}\label{lem:triviality-01}Let $T\in\mathcal{L}(H)$ have the polar decomposition $T=U_T|T|$, and let $A\in\mathcal{L}(H)$ be given such that $\big[A, |T|\big]=0$. Then $[A,T]=0$ if and only if  $[A,U_T]=0$.
\end{lemma}
\begin{proof}
Suppose that $AU_T=U_T A$. Then
$$AT=AU_T |T|=U_T A|T|=U_T|T|A=TA.$$
Conversely, assume that $AT=TA$. Then by Lemma~\ref{lem:rang characterization-1}
\begin{align*}\overline{\mathcal{R}(AU_T)}=\overline{\mathcal{R}(AT)}=\overline{\mathcal{R}(TA)}\subseteq \overline{\mathcal{R}(T)},
\end{align*}
which means that
\begin{equation}\label{equ:range inclusion-001}U_TU_T^* AU_T=AU_T.
\end{equation}
Since $A|T|=|T|A$, we see from Lemma~\ref{lem:commutative property extended-3}~(ii) that
\begin{equation}\label{equ:commutatity of A-001}AU_T^*U_T=U_T^*U_TA.
\end{equation}
Furthermore, due to $[A,T]=\big[A,|T|\big]=0$ we have
\begin{align*}&U_T^*AU_T|T|=U_T^*AT=U_T^*TA=|T|A=A|T|,
\end{align*}
so
$U_T^*AU_Tx=Ax$
for every $x\in\overline{\mathcal{R}(|T|)}=\mathcal{R}(U_T^*)$. Hence,
$U_T^*AU_TU_T^*=AU_T^*$. This, together with \eqref{equ:commutatity of A-001}, yields
$$U_T^*AU_T=(U_T^*AU_TU_T^*)U_T=AU_T^*U_T=U_T^*U_T A,$$
which is combined with \eqref{equ:range inclusion-001} to get
\begin{align*}AU_T=U_T(U_T^*AU_T)=U_T(U_T^*U_T A)=U_TA. \qquad \qedhere
\end{align*}
\end{proof}

An application of  Theorem~\ref{thm:polar decom of three operators} can now be provided  as follows.
\begin{theorem}\label{thm:application of main thm} Let  $T=U_T|T|$, $S=U_S|S|$ and $A=U_A|A|$ be the polar decompositions of $T,A,S\in\mathcal{L}(H)$, respectively. Suppose that  $\big[|T|,A\big]=0$. Let $X$ and $Y$ be defined by \eqref{equ:defn of X and Y}, and let $W=U_T^*U_T U_AU_SU_S^*$. Then the following statements are  equivalent:
\begin{enumerate}
\item[{\rm (i)}] $X=U_TU_AU_S|X|$  is the polar decomposition;
\item[{\rm (ii)}]$Y=W|Y|$  is the polar decomposition;
\item[{\rm (iii)}] $X=U_TU_AU_S|X|$;
\item[{\rm (iv)}]$Y=W|Y|$;
\item[{\rm (v)}] $|T|\cdot |A|\cdot |S^*|$ is positive.
\end{enumerate}
\end{theorem}
\begin{proof}  Put $Z=|T|\cdot |A|\cdot |S^*|$.
By $\big[|T|,A\big]=0$, we have $\big[|T|,A^*\big]=0$, which in turn gives $\big[|T|,A^*A\big]=0$, hence $\big[|T|,|A|\big]=0$ by  Lemma~\ref{lem:commutative property extended-3}~(i). Therefore,  by the definitions of $Y$ and $Z$ we see that $|Y|=|Z|$. Furthermore,
due to $\big[|T|,|A|\big]=\big[|T|,A\big]=0$ we have $\big[|T|,U_A\big]=0$ (see Lemma~\ref{lem:triviality-01}), and thus $\big[|T|,U_A^*\big]=0$.
So  by Lemma~\ref{lem:commutative property extended-3}~(ii) we have
\begin{equation*}[U_T^*U_T, U_A]=[U_T^*U_T, U_A^*]=0.
\end{equation*}
It follows that
\begin{equation}\label{equ:two equation-needed-01}U_A^*Y=|T|\cdot U_A^* A\cdot |S^*|=Z, \quad U_AU_A^*Y=|T|\cdot U_AU_A^* A\cdot |S^*|=Y,\end{equation}
and
\begin{align*}\overline{\mathcal{R}(|T|\cdot |A|)}=&\overline{\mathcal{R}(|T|U_A^*U_A)}=\overline{\mathcal{R}(U_A^*U_A|T|)}=\overline{\mathcal{R}(U_A^*U_AU_T^*U_T)}\\
=&\overline{\mathcal{R}(U_A^*U_T^*U_TU_A)}=\overline{\mathcal{R}\big[(U_T U_A)^*(U_T U_A)\big]}=\overline{\mathcal{R}\big[(U_T U_A)^*\big]}.
\end{align*}
Thus,
\begin{equation}\label{equ:range closure-needed-01}\overline{\mathcal{R}(|T|\cdot |A|)}=\overline{\mathcal{R}(U_A^*U_T^*)}=\overline{\mathcal{R}(U_A^*U_T^*U_T)},
\end{equation}
since it is clear that $\mathcal{R}(U_T^*)=\mathcal{R}(U_T^*U_T)$.

From Theorem~\ref{thm:polar decom of three operators}  and Lemma~\ref{lem:technical lemma for three operators} we see that
(i)$\Longleftrightarrow$(ii) and (iii)$\Longleftrightarrow$(iv).
The implication (iv)$\Longrightarrow$(v) can be derived directly from (ii)$\Longrightarrow$(iii) in Lemma~\ref{lem:technical lemma for three operators} and the first two equations in \eqref{equ:two equation-needed-01}.

(v)$\Longrightarrow$(ii). Assume that $Z\ge 0$. As is shown above, we have $|Y|=|Z|=Z$. In view of \eqref{equ:two equation-needed-01} and Theorem~\ref{thm:polar decom of three operators}, it needs only to prove that $\overline{\mathcal{R}(Z)}=\overline{\mathcal{R}(W^*)}$.
Actually, from $Z=Z^*$ and $[|T|,|A|]=0$ we can obtain
$$(|T|\cdot |A|)|S^*|=|S^*|(|T|\cdot |A|),$$
which leads by Lemma~\ref{lem:commutative property extended-3}~(ii) to
$$(|T|\cdot |A|)U_SU_S^*=U_SU_S^*(|T|\cdot |A|).$$
It follows from Lemma~\ref{lem:rang characterization-1} and \eqref{equ:range closure-needed-01} that
\begin{align*}\overline{\mathcal{R}(Z)}=&\overline{\mathcal{R}(|T|\cdot |A|U_SU_S^*)}=\overline{\mathcal{R}(U_SU_S^*|T|\cdot |A|)}\\
=&\overline{\mathcal{R}(U_SU_S^*U_A^*U_T^*U_T)}=\overline{\mathcal{R}(W^*)}.\qedhere
\end{align*}

\end{proof}

\begin{remark} In the special case that $A=I_H$, the equivalence of (i), (iii) and (v) in the preceding theorem
was given in \cite[Theorem~3.5]{LLX-BJMA}.
\end{remark}

\begin{definition}\label{rem:some results of widetilde U} Suppose that $T,A\in\mathcal{L}(H)$  have the polar decompositions $T=U_T|T|$ and $A=U_A|A|$. For every $\alpha>0, \beta>0$, let
\begin{equation}\label{equ:defn of widetilde U}X^{\alpha, \beta}_{T,A}=|T|^{\alpha}\cdot  A\cdot U_T|T|^{\beta},\ Y^{\alpha, \beta}_{T,A}=|T|^{\alpha}\cdot  A\cdot |T^*|^\beta, \ W_{T,A}=U_T^*U_TU_AU_T.
\end{equation}
\end{definition}

\begin{corollary}\label{cor:equivalent conditions of binormal} Let $T,A\in\mathcal{L}(H)$ have the polar decompositions
$T=U_T|T|$ and $A=U_A|A|$. Suppose that  $\big[|T|,A\big]=0$. Let $X^{\alpha, \beta}_{T,A}$, $Y^{\alpha, \beta}_{T,A}$ and $W_{T,A}$ be defined by \eqref{equ:defn of widetilde U} for $\alpha>0$ and $\beta>0$. Then the following statements are equivalent:
\begin{enumerate}
\item[{\rm (i)}]$X^{\alpha, \beta}_{T,A}=W_{T,A}\left|X^{\alpha, \beta}_{T,A}\right|$ is the polar decomposition;
\item[{\rm (ii)}]$Y^{\alpha, \beta}_{T,A}=W_{T,A}U_T^*\left|Y^{\alpha, \beta}_{T,A}\right|$ is the polar decomposition;
\item[{\rm (iii)}]$X^{\alpha, \beta}_{T,A}=W_{T,A}\left|X^{\alpha, \beta}_{T,A}\right|$;
\item[{\rm (iv)}]$Y^{\alpha, \beta}_{T,A}=W_{T,A}U_T^*\left|Y^{\alpha, \beta}_{T,A}\right|$;
\item[{\rm (v)}] $|T|^\alpha \cdot |A|\cdot |T^*|^\beta$ is positive.
\end{enumerate}
If in addition $\big[|T^*|,|A|\big]=0$, and one of (i)--(v) is satisfied, then (i)--(v) are all valid when $\alpha$ and $\beta$ are replaced by any positive numbers $\alpha^\prime$ and $\beta^\prime$, respectively.
\end{corollary}
\begin{proof}(1) Put $\widetilde{T}=|T|^{\alpha}$ and $S=U_T|T|^{\beta}$.
By Lemma~\ref{lem:Range Closure of T alpha and T}, $\widetilde{T}=U_T^*U_T\big|\widetilde{T}\big|$ is the polar decomposition. It follows from the proof of \cite[Theorem~4.14]{LLX-BJMA} that $S=U_T|T|^{\beta}$ is the polar decomposition such that $|S^*|=|T^*|^\beta$. Thus,
\begin{equation*}|\widetilde{T}|=|T|^\alpha,\ |S^*|=|T^*|^\beta,\ U_{\widetilde{T}}U_AU_S=W_{T,A},\ U_{\widetilde{T}}^* U_{\widetilde{T}}U_AU_SU_S^*=W_{T,A}U_T^*.\end{equation*}
Since $\big[|T|,|A|\big]=0$, by Lemma~\ref{lem:commutative property extended-3}~(i) we have $\big[|\widetilde{T}|,A\big]= \big[|T|^\alpha, A\big]=0$.
Hence, the equivalence of (i)--(v) follows from Theorem~\ref{thm:application of main thm}.

(2) Suppose furthermore that $\big[|T^*|,|A|\big]=0$. Then for every $\alpha^\prime>0$ and $\beta^\prime>0$, by Lemma~\ref{lem:commutative property extended-3}~(i) we have
$\big[|T|^{\alpha^\prime},|A|\big]=\big[|T^*|^{\beta^\prime},|A|\big]=0.$
Therefore, both $|T|^{\alpha^\prime}\cdot |A|$ and $|A|\cdot |T^*|$ are positive (see Lemma~\ref{lem:commutative property extended-0}).
Combining Lemma~\ref{lem:commutative property extended-0} and Lemma~\ref{lem:commutative property extended-3}~(i), we see that
\begin{align*}|T|^{\alpha^\prime} \cdot |A|\cdot |T^*|^{\beta^\prime}\ge 0&\Longleftrightarrow \left[|T|^{\alpha^\prime} \cdot |A|,|T^*|^{\beta^\prime}\right]=0\Longleftrightarrow \left[|T|^{\alpha^\prime} \cdot |A|,|T^*|\right]=0\\
&\Longleftrightarrow |T|^{\alpha^\prime} \cdot |A|\cdot |T^*|\ge 0\Longleftrightarrow \left[|A|\cdot |T^*|, |T|^{\alpha^\prime}\right]=0\\
&\Longleftrightarrow \big[|A|\cdot |T^*|, |T|\big]=0\Longleftrightarrow |T|\cdot |A|\cdot |T^*|\geq 0.
\end{align*}
Consequently,
$$|T|^{\alpha^\prime} \cdot |A|\cdot |T^*|^{\beta^\prime}\ge 0\Longleftrightarrow |T|^{\alpha}\cdot |A|\cdot |T^*|^{\beta}\ge 0.\qedhere$$
\end{proof}
\begin{remark}In the  special case that $A=I_H$, the equivalence of (i), (iii) and (v) in the preceding corollary can be found in
\cite[Theorem~4.14]{LLX-BJMA}.
\end{remark}

\begin{remark}There exist non-zero matrices $T,A\in\mathbb{C}^3$ such that $A\ne I_3$, $\big[|T|,A\big]=\big[|T^*|,|A|\big]=0$, whereas $\big[|T|,|T^*|\big]\ne 0$. Such an example is as follow.
\end{remark}

\begin{example} Let $T,A\in\mathbb{C}^3$ be given by
\begin{align*} T = \left(
                                                                    \begin{array}{ccc}
                                                                      1 & 0 & 1 \\
                                                                      0 & 0 & 0 \\
                                                                      0 & 0 & 1 \\
                                                                    \end{array}
                                                                  \right)\quad\mbox{and}\quad
A = \left(
                                                                    \begin{array}{ccc}
                                                                      1 & 0 & 0 \\
                                                                      0 & 2 & 0 \\
                                                                      0 & 0 & 1 \\
                                                                    \end{array}
                                                                  \right).
\end{align*}
Since $A=A^*$ and $[A,T]=0$, we have $\big[|T|,A\big]=\big[|T^*|,|A|\big]=0$.  A simple calculation shows that
\begin{align*}
|T|=\left(
                                                                    \begin{array}{ccc}
                                                                      \frac{2\sqrt{5}}{5} & 0 & \frac{\sqrt{5}}{5} \\
                                                                      0 & 0 & 0 \\
                                                                      \frac{\sqrt{5}}{5} & 0 & \frac{3\sqrt{5}}{5} \\
                                                                    \end{array}
                                                                  \right) \quad\mbox{and}\quad
|T^*|=\left(
                                                                    \begin{array}{ccc}
                                                                      \frac{3\sqrt{5}}{5} & 0 & \frac{\sqrt{5}}{5} \\
                                                                      0 & 0 & 0 \\
                                                                      \frac{\sqrt{5}}{5} & 0 & \frac{2\sqrt{5}}{5} \\
                                                                    \end{array}
                                                                  \right),
\end{align*}
which gives $\big[|T|,|T^*|\big]\ne 0$.
\end{example}

\vspace{2ex}

\end{document}